\documentclass[oneside,reqno]{amsart}
\usepackage{amssymb}
\usepackage{amscd}
\usepackage{amsmath}
\usepackage{array}
\usepackage{longtable}
\usepackage{subcaption}
\usepackage{graphicx}
\usepackage{color,soul}
\usepackage{multicol}
\numberwithin{equation}{section}
\newtheorem{thm}{Theorem}
\newtheorem{pro}{Property}
\newtheorem{defn}{Definition}
\newtheorem{exa}{Example}

\begin{document}

\title[NEW INTEGRAL TRANSFORM FOR SOLVING DIFFERENTIAL EQUATIONS]{NEW INTEGRAL TRANSFORM: SHEHU TRANSFORM A GENERALIZATION OF SUMUDU AND LAPLACE TRANSFORM FOR SOLVING DIFFERENTIAL EQUATIONS}



\author[]{\bfseries SHEHU MAITAMA$^{*}$, WEIDONG ZHAO}

\address{School of Mathematics, Shandong University, Jinan, Shandong 250100, P.R. China}
\address{$^*$Corresponding author: \textnormal{smusman12@sci.just.edu.jo}}

\begin{abstract}
In this paper, we introduce a Laplace-type integral transform
called the Shehu transform which is a generalization of the
Laplace and the Sumudu integral transforms for solving
differential equations in the time domain. The proposed integral
transform is successfully derived from the classical Fourier
integral transform and is applied to both ordinary and partial
differential equations to show its simplicity, efficiency, and the
high accuracy.
\end{abstract}
\maketitle

\section{Introduction}

Historically, the origin of the integral transforms can be traced
back to the work of P. S. Laplace in 1780s and Joseph Fourier in
1822. In recent years, differential and integral equations have
been solved using many integral transforms (\cite{R1}-\cite{R11}).
The Laplace transform, and Fourier integral transforms are the
most commonly used in the literature. The Fourier integral
transform \cite{R12} was named after the French mathematician
Joseph Fourier. Mathematically, Fourier integral transform is
defined as:
\begin{equation}\label{eq1.1}
\digamma[f(t)]=f(\omega)=\frac{1}{\sqrt{2\pi}}\int_{-\infty}^{\infty}\exp\left(-i\omega
t\right)f(t)dt.
\end{equation}
The Fourier transform have many applications in physics and
engineering processes \cite{R13}. The Laplace integral transform
is similar with the Fourier transform and is defined as:
\begin{equation}
\pounds[f(t)]=F(s)=\int_{-\infty}^{\infty}\exp\left(-st\right)f(t)dt.
\end{equation}
The Laplace transform is highly efficient for solving some class
of ordinary and partial differential equations \cite{R14}. By
replacing the variable $i\omega$ with the variable $s$ in
Equ.(\ref{eq1.1}), the well-known Fourier transform will become a
Laplace transform and the vice-versa. The only difference between
the Laplace transform, and the Fourier transform is that the
Laplace transform can be defined for both stable and unstable
system while the Fourier transform can only be defined on a stable
system. In mathematical literature, the discrete-time equivalent
of the Laplace transform called z-transform \cite{R15} converts a
discrete-time signal into a complex frequency-domain
representation. The basic idea of the z-transform was known to
Laplace and later it was re-introduced by the Jewish-Polish
mathematician Witold Hurewicz to treat a sampled-data control
systems used with radar in 1947 (\cite{R16}-\cite{R17}). In
mathematics and signal processing, the bilateral or two-sided
z-transform of a discrete-time signal $x[n]$ is the normal power
series $X(z)$ which is defined as:
\begin{equation}
X(z)=Z\{x[n]\}=\sum_{n=-\infty}^{\infty}x[n]z^{-n},
\end{equation}
where $n$ is an integer and $z$ is in general a complex number
\cite{R18}.

The multiplicative version of the two-sided Laplace transform
called the Mellin integral transform is defined as \cite{R19}:
\begin{equation}
M[f(s);s]=f^{*}(s)=\int_{0}^{\infty}x^{s-1}f(x)dx.
\end{equation}

The Mellin integral transform is similar with the Laplace
transform and Fourier transform and is widely applied in computer
science and number theory due to its invariant property
\cite{R20,R21}. In railway engineering, the Laplace-Carson
transform \cite{R22} which is a Laplace-type integral transform
named after Pierre Simon Laplace and John Renshaw Carson is
defined as:

\begin{equation}
\hat{f}_C(p)=p\int_{0}^{\infty}\exp\left(-pt\right)f(t)dt,\,\,t\geq0.
\end{equation}
The Laplace-Carson integral transform have many applications in
physics and engineering and can easily  be converted into a Mellin
deconvolution problem, see (\cite{R23},\cite{R24}). In
mathematics, the Hankel's integral transform \cite{R25} which is
similar to the Fourier transform was first introduced by the
German mathematician Hermann Hankel and was widely used in
physical science and engineering \cite{R26}. The Hankel's
transform is defined as:
\begin{equation}
F_v(s)=H_v[f(r)]=\int_{0}^{\infty}rf(r)J_v(sr)dr,\,\,\,\,\,r\geq0,
\end{equation}
where $J_v$ is the Bessel function of the first kind of order $v$
with $v\geq-\frac{1}{2}$.

In 1993, Watugala introduced a Laplace-like integral transform
called the Sumudu integral transform \cite{R27}. In recent years,
Sumudu transform has been applied to many real-life problems
because of its scale and unit preserving properties
(\cite{R28}-\cite{R31}). The mathematical definition of the Sumudu
transform is given by:
\begin{equation}
\textbf{S}[f(t)](u)=G(u)=\frac{1}{u}\int_{0}^{\infty}\exp\left(\frac{-t}{u}\right)f(t)dt,
\end{equation}
provided the integral exists for some $u$. Based on the basic idea
of the Laplace and the Sumudu integral transform, the Elzaki
transform was proposed in 2011. The Elzaki transform is closely
related with the Laplace transform, Sumudu transform, and the
natural transform. Elzaki transform is defined as \cite{R32}:
\begin{equation}
E[f(t)]=T(u)=u\int_{0}^{\infty}\exp\left(\frac{-t}{u}\right)f(t)dt,
\end{equation}
provided the integral exists for some $u$.

The natural transform \cite{R33} which is similar to Laplace and
Sumudu integral transform was introduced in 2008. In recent years,
natural transform was successfully applied to many applications
(see \cite{R34,R35}). The natural transform is defined by the
following integral:
\begin{equation}
{\mathbb
N^+}[f(t)](s,u)=R(s,u)=\frac{1}{u}\int_{0}^{\infty}\exp\left(\frac{-st}{u}\right)f(t)dt,
\, \, s\,>0,\,\, u>0,
\end{equation}
provided the integral exists for some variables $u$ and $s$.
Recently, a new integral transform called the ${\mathbb
M}$-transform which is also similar to natural transform is
introduced by Srivastava et al. in 2015. Mathematically speaking,
${\mathbb M}$-transform is closely connected with the well-known
Laplace transform and the Sumudu integral transform. ${\mathbb
M}$-transform was successfully applied to first order
initial-boundary value problem (see Srivastava et al. \cite{R36}).
The ${\mathbb M}$-transform is defined as:
\begin{equation}
{\mathbb
M_{\rho,m}}[f(t)](u,v)=\int_{0}^{\infty}\frac{\exp\left(-ut\right)f(vt)}{\left(t^m+v^m\right)^\rho}dt,
\end{equation}
$\left(\rho \in {\mathbb C};\,\Re(\rho)\geq0,\,m\in {\mathbb
Z_+}={1,2,3,\cdots}\right),$ where both $u\in {\mathbb C}$ and
$v\in {\mathbb R_+}$ are the ${\mathbb M}$-transform variables.

In 2013, Atangana and Kilicman introduced a novel integral
transform called the Abdon-Kilicman integral transform \cite{R37}
for solving some differential equations with some kind of
singularities. The novel integral transform is defined as:
\begin{equation}
M_n(s)=M_n[f(x)](s)=\int_{0}^{\infty}x^n\exp\left(-xs\right)f(x)dx.
\end{equation}
The Atangana-Kilicman integral becomes Laplace transform when
$n=0$. Recently, a Laplace-type integral transform called the Yang
transform (\cite{R38}-\cite{R40}) for solving steady heat transfer
problems was introduced in 2016. The integral transform is defined
as:
\begin{equation}
Y[\phi(\tau)]=\phi(\omega)=\int_{0}^{\infty}\exp\left(\frac{-\tau}{\omega}\right)\phi(\tau)d\tau,
\end{equation}
provided the integral exists for some $\omega$.

Due to the rapid development in the physical science and
engineering models, there are many other integral transforms in
the literature. However, most of the existing integral transforms
have some limitations and cannot be used directly to solved
nonlinear problems or many complex mathematical models. As a
result, many authors became highly interested to come up with the
alternative approach for solving many real-life problems. In 2016,
Atangana and Alkaltani introduced a new double integral equation
and their properties based on the Laplace transform and
decomposition method. The double integral transform was
successfully applied to second order partial differential equation
with singularity called the two-dimensional Mboctara equation
\cite{R41}. Recently, Eltayeb applied double Laplace decomposition
method to nonlinear partial differential equations \cite{R42}. In
2017, Belgacem el at. extended the applications of the natural and
the Sumudu transforms to fractional diffusion and Stokes fluid
flow realms \cite{R43}.

Motivated by the above-mentioned researches, in this paper we
proposed a Laplace-type integral transform called Shehu transform
for solving both ordinary and partial differential equations. The
Laplace-type integral transform converges to Laplace transform
when $u=1$, and to Yang integral transform when $s=1$. The
proposed integral transform is successfully applied to both
ordinary and partial differential equations. All the results
obtained in the applications section can easily be verified using
the Laplace or Fourier integral transforms. Throughout this paper,
the Shehu transform is denoted by an operator ${\mathbb S[.]}$.

\vspace{1cm}
\section{Main result}
\begin{defn}\label{D1}
The Shehu transform of the function $v(t)$ of exponential order is defined over the set of functions,\\
$A=\left\{v(t):\exists \, \, N,\, \eta _{1} ,\, \eta _{2}
>0,\, \, \left|v(t)\right|<N\exp\left(\frac{\left|t\right|}{\eta _{i} }\right)
,\, \, \textnormal{if}\, \, t\in (-1)^{i} \times \left[0,\infty
\right)\right\}$,\\by the following integral
\begin{eqnarray}\label{eq2.1}
{\mathbb S}\left[v(t)\right]&=&V(s,u)=\int _{0 }^{\infty
}\exp\left(\frac{-st}{u}\right)v(t)dt\nonumber\\&=&\lim_{\alpha\rightarrow\infty}\int
_{0 }^{\alpha }\exp\left(\frac{-st}{u}\right)v(t)dt;\, \, \, \,
s\,>0,\,\, u>0.
\end{eqnarray}
\end{defn}
It converges if the limit of the integral exists, and diverges if
not.

The inverse Shehu transform is given by
\begin{equation}
{\mathbb S^{-1}}\left[V(s,u)\right]=v(t),\, \, \,for\,\,\, t\geq0.
\end{equation}
Equivalently

\begin{equation}\label{eq2.3}
v(t)={\mathbb S^{-1}}\left[V(s,u)\right]=\frac{1}{2\pi i}\int
_{\alpha-i\infty}^{\alpha+i\infty}\frac{1}{u}\exp\left(\frac{st}{u}\right)V(s,u)\,
ds,
\end{equation}
where $s$ and $u$ are the Shehu transform variables, and $\alpha$
is a real constant and the integral in Equ.(\ref{eq2.3}) is taken
along $s=\alpha$ in the complex plane $s=x+iy$.

\begin{thm}
\textbf{The sufficient condition for the existence of Shehu
transform}.
 If the function $v(t)$ is piecewise continues in every finite
interval $0\leq t\leq\beta$ and of exponential order $\alpha$ for
$t>\beta$. Then its Shehu transform $V(s,u)$ exists.
\end{thm}
\textbf{Proof.} For any positive number $\beta$, we algebraically
deduce
\begin{equation}
\int_{0}^{\infty}\exp\left(-\frac{st}{u}\right)v(t)dt=\int_{0}^{\beta}\exp\left(-\frac{st}{u}\right)v(t)dt+\int_{\beta}^{\infty}\exp\left(-\frac{st}{u}\right)v(t)dt.
\end{equation}
Since the function $v(t)$ is piecewise continues in every finite
interval $0\leq t\leq\beta$, then the first integral on the right
hand side exists. Besides, the second integral on the right hand
side exists, since the function $v(t)$ is of exponential order
$\alpha$ for $t>\beta$. To verify this claim, we consider the
following case
\begin{eqnarray*}
\left|\int_{\beta}^{\infty}\exp\left(-\frac{st}{u}\right)v(t)dt\right|&\leq&\int_{\beta}^{\infty}\left|\exp\left(-\frac{st}{u}\right)v(t)\right|dt
\\&\leq&\int_{0}^{\infty}\exp\left(-\frac{st}{u}\right)\left|v(t)\right|dt\\&\leq&\int_{\beta}^{\infty}\exp\left(-\frac{st}{u}\right)N\exp\left(\alpha t\right)dt
\\&=&N\int_{\beta}^{\infty}\exp\left(-\frac{(s-\alpha
u)t}{u}\right)dt\\&=&-\frac{u N}{(s-\alpha
u)}\lim_{\gamma\rightarrow\infty}\left[\exp\left(-\frac{(s-\alpha
u)t}{u}\right)\right]^\gamma_0\\&=&\frac{u N}{s-\alpha u}.
\end{eqnarray*}
The proof is complete.\qquad\qquad $\Box$

 \begin{pro}\label{P1}
 \textbf{Linearity property of Shehu transform}. Let the functions $\alpha v(t)$ and $\beta w(t)$ be in set A, then  $\left(\alpha v(t)+\beta w(t)\right)\in $A, where $\alpha$ and $\beta$ are nonzero arbitrary constants, and
 \begin{equation}
 {\mathbb S}\left[\alpha v(t)+\beta w(t)\right]=\alpha{\mathbb S}\left[v(t)\right]+\beta{\mathbb
 S}\left[w(t)\right].
 \end{equation}
\textbf{Proof.} Using the Definition 1 of Shehu transform, we get
\begin{eqnarray*}
{\mathbb S}\left[\alpha v(t)+\beta w(t)\right]&=&\int _{0
}^{\infty }\exp\left(\frac{-st}{u}\right)(\alpha v(t)+\beta w(t))dt\\
&=&\int _{0 }^{\infty }\exp\left(\frac{-st}{u}\right)(\alpha v(t))
\, dt+\int
_{0 }^{\infty }\exp\left(\frac{-st}{u}\right)(\beta w(t))dt\\
&=&\alpha\int _{0 }^{\infty
}\exp\left(\frac{-st}{u}\right)v(t)dt+\beta\int _{0 }^{\infty
}\exp\left(\frac{-st}{u}\right)w(t)dt\\
&=&\alpha u\int _{0 }^{\infty }\exp\left(-st\right)v(ut)dt+\beta
u\int _{0 }^{\infty }\exp\left(-st\right)
w(ut)dt\\
&=&\alpha {\mathbb S}\left[v(t)\right]+\beta {\mathbb
S}\left[w(t)\right].
\end{eqnarray*}
The proof is complete.\qquad\qquad $\Box$
\end{pro}
In particular, using the Definition \ref{D1} and Property
\ref{P1}, we obtain
\begin{equation}
 {\mathbb S}\left[3\cos(t)+5\sin(2t)\right]= 3{\mathbb S}\left[\cos(t)\right]+ 5{\mathbb
 S}\left[\sin(2t)\right]\nonumber
 \end{equation}
 \begin{equation}
 =\frac{3us}{s^2+u^2}+\frac{5u^2}{s^2+(2u)^2}\nonumber,
 \end{equation}
see entries of table 1.

 \begin{pro}\label{P2}
 \textbf{Change of scale property of Shehu transform.} Let the function $v(\beta
 t)$ be in set A, where $\beta$ is an arbitrary constant. Then
 \begin{equation}
 {\mathbb S}\left[v(\beta
 t)\right]=\frac{u}{\beta}V\left(\frac{s}{\beta},u\right).
 \end{equation}
\textbf{Proof.} Using the Definition 1 of Shehu transform, we
deduce
\begin{eqnarray}\label{eq2.7}
{\mathbb S}\left[v(\beta t)\right]&=&\int _{0 }^{\infty
}\exp\left(\frac{-st}{u}\right) v(\beta t)dt
\end{eqnarray}
Substituting $\eta=\beta t$ which implies $t=\frac{\eta}{\beta}$
and $dt=\frac{d\eta}{\beta}$ in Equ.(\ref{eq2.7}) yields
\begin{eqnarray*}
{\mathbb S}\left[v(\beta t)\right]&=&\frac{1}{\beta}\int _{0
}^{\infty }\exp\left(\frac{-s\eta}{u\beta}\right)
v(\eta)d\eta\\
&=&\frac{1}{\beta}\int _{0 }^{\infty
}\exp\left(\frac{-st}{u\beta}\right)
v\left(t\right)dt\\
&=&\frac{u}{\beta}\int _{0 }^{\infty
}\exp\left(\frac{-st}{\beta}\right)
v(ut)dt\\
&=&\frac{u}{\beta}V\left(\frac{s}{\beta},u\right).
\end{eqnarray*}
The proof is complete.\qquad\qquad $\Box$
\end{pro}
\begin{thm}\textbf{Derivative of Shehu transform.}
 If the function $v^{(n)}(t)$ is
the $nth$ derivative of the function $v(t)\in A$ with respect to
$'t'$, then its Shehu transform is defined by
\begin{equation}\label{eq2.8}
{\mathbb S}
\left[v^{(n)}(t)\right]=\frac{s^n}{u^n}V(s,u)-\sum_{k=0}^{n-1}\left(\frac{s}{u}\right)^{n-(k+1)}v^{(k)}(0).
\end{equation}
\end{thm}
When n=1,\,\,2,\,\,and 3 in Equ.\;(\ref{eq2.8}) above, we obtain
the following derivatives with respect to $t$.
\begin{equation}\label{eq2.9}
{\mathbb S} \left[v'(t)\right]=\frac{s}{u} V(s,u)-v(0).
\end{equation}
\begin{equation} {\mathbb
S}\left[v''(t)\right]=\frac{s^{2} }{u^{2} }
V(s,u)-\frac{s}{u}v(0)-v'(0).
\end{equation}
\begin{equation} {\mathbb
S}\left[v'''(t)\right]=\frac{s^{3} }{u^{3} } V(s,u)-\frac{s^2}{u^2
}v(0)-\frac{s}{u}v'(0)-v''(0).
\end{equation}
\textbf{Proof.} Now suppose Equ.\;(\ref{eq2.8}) is true for $n=k$.
Then using Equ.\;(\ref{eq2.9}) and the induction hypothesis, we
deduce
\begin{eqnarray*}\label{eq12}
{\mathbb S} \left[(v^{(k)}(t))'\right]&=&\frac{s}{u}
{\mathbb S} \left[v^{(k)}(t)\right]-v^{(k)}(0)\\
&=&\frac{s}{u}\left[\frac{s^k}{u^k}{\mathbb S}\left[v(t)\right]
-\sum_{i=0}^{k-1}\left(\frac{s}{u}\right)^{k-(i+1)}v^{(i)}(0)\right]-v^{(k)}(0)\\
&=&\left(\frac{s}{u}\right)^{k+1}{\mathbb S}
\left[v(t)\right]-\sum_{i=0}^{k}\left(\frac{s}{u}\right)^{k-i}v^{(i)}(0),
\end{eqnarray*}
which implies that Equ. (\ref{eq2.8}) holds for $n=k+1$. By
induction hypothesis the proof is complete. \qquad\qquad $\Box$\\
The following important properties are obtain using the Leibniz's
rule
\begin{eqnarray*}
{\mathbb S}\left[\frac{\partial v(x,t)}{\partial x}\right]=\int
_{0 }^{\infty }\exp\left(\frac{-st}{u}\right)\frac{\partial
v(x,t)}{\partial x}\, dt=\frac{\partial}{\partial x}\int _{0
}^{\infty }\exp\left(\frac{-st}{u}\right)v(x,t)\, dt\\
=\frac{\partial}{\partial x}\left[V(x,s,u)\right]\Rightarrow
{\mathbb S}\left[\frac{\partial v(x,t)}{\partial
x}\right]=\frac{d}{d x}\left[V(x,s,u)\right],
\end{eqnarray*}
\begin{eqnarray*}
{\mathbb S}\left[\frac{\partial^2 v(x,t)}{\partial
x^2}\right]=\int _{0 }^{\infty
}\exp\left(\frac{-st}{u}\right)\frac{\partial^2 v(x,t)}{\partial
x^2}\, dt=\frac{\partial^2}{\partial x^2}\int _{0
}^{\infty }\exp\left(\frac{-st}{u}\right)v(x,t)\, dt\\
=\frac{\partial^2}{\partial x^2}\left[V(x,s,u)\right]\Rightarrow
{\mathbb S}\left[\frac{\partial^2 v(x,t)}{\partial
x^2}\right]=\frac{d^2}{d x^2}\left[V(x,s,u)\right],
\end{eqnarray*}
and
\begin{eqnarray*}
{\mathbb S}\left[\frac{\partial^n v(x,t)}{\partial
x^n}\right]=\int _{0 }^{\infty
}\exp\left(\frac{-st}{u}\right)\frac{\partial^n v(x,t)}{\partial
x^n}\, dt=\frac{\partial^n}{\partial x^n}\int _{0
}^{\infty }\exp\left(\frac{-st}{u}\right)v(x,t)dt\\
=\frac{\partial^n}{\partial x^n}\left[V(x,s,u)\right]\Rightarrow
{\mathbb S}\left[\frac{\partial^n v(x,t)}{\partial
x^n}\right]=\frac{d^n}{d x^n}\left[V(x,s,u)\right].
\end{eqnarray*}
\vskip 2mm
\section{Some useful results of Shehu transform}
\vskip 2mm
\begin{pro}
Let the function $v(t)=1$ be in set A. Then its Shehu transform is
given by
\begin{equation}
{\mathbb S}\left[1\right]=\frac{u}{s}.
\end{equation}
\textbf{Proof.} Using Equ.(\ref{eq2.1}), we deduce
\begin{equation}
{\mathbb
S}\left[1\right]=\int_{0}^{\infty}\exp\left(\frac{-st}{u}\right)dt=-\frac{u}{s}\lim_{\gamma\rightarrow\infty}\left[\exp\left(\frac{-st}{u}\right)\right]^\gamma_0=\frac{u}{s}\nonumber.
\end{equation}
This ends the proof.\qquad\qquad $\Box$
\end{pro}
\begin{pro}
Let the function $v(t)=t$ be in set A. Then its Shehu transform is
given by
\begin{equation}
{\mathbb S}\left[t\right]=\frac{u^2}{s^2}.
\end{equation}
\textbf{Proof.} Using the Definition 1 of the Shehu transform and
integration by parts, we get
\begin{eqnarray*}
{\mathbb
S}\left[t\right]&=&\int_{0}^{\infty}t\exp\left(\frac{-st}{u}\right)dt=\frac{u}{s}\lim_{\gamma\rightarrow\infty}\left[t\exp\left(\frac{-st}{u}\right)\right]^\gamma_0+\frac{u}{s}\int_{0}^{\infty}\exp\left(\frac{-st}{u}\right)dt\\&&=-\frac{u^2}{s^2}\lim_{\gamma\rightarrow\infty}\left[\exp\left(\frac{-st}{u}\right)\right]^\gamma_0=\frac{u^2}{s^2}.
\end{eqnarray*}
Thus the proof ends.\qquad\qquad $\Box$
\end{pro}
\begin{pro}
Let the function $v(t)=\frac{t^{n}}{n!}\,\,\,n=0,1,2..$ be in set
A. Then its Shehu transform is given by
\begin{equation}
{\mathbb
S}\left[\frac{t^{n}}{n!}\right]=\left(\frac{u}{s}\right)^{n+1}.
\end{equation}
\vskip 2mm \textbf{Proof.} From the Definition 1 of the Shehu
transform and integration by parts, we deduce
\begin{eqnarray*}
&&{\mathbb S}\left[t^{n}\right]=\int_{0}^{\infty}t^n\exp\left(\frac{-st}{u}\right)dt=\frac{u}{s}n\int_{0}^{\infty}t^{n-1}\exp\left(\frac{-st}{u}\right)dt\\&&=\frac{u^2}{s^2}n(n-1)\int_{0}^{\infty}t^{n-2}\exp\left(\frac{-st}{u}\right)dt\\
&=&\frac{u^3}{s^3}n(n-1)(n-2)\int_{0}^{\infty}t^{n-3}\exp\left(\frac{-st}{u}\right)dt\\&&=\frac{u^4}{s^4}n(n-1)(n-2)(n-3)\int_{0}^{\infty}t^{n-4}\exp\left(\frac{-st}{u}\right)dt\\
&=&\frac{u^5}{s^5}n(n-1)(n-2)(n-3)(n-4)\int_{0}^{\infty}t^{n-5}\exp\left(\frac{-st}{u}\right)dt=
\cdots = n!\left(\frac{u}{s}\right)^{n+1}.
\end{eqnarray*}
The proof is completed.\qquad\qquad $\Box$
\end{pro}
\begin{pro}
Let the function
$v(t)=\frac{t^{n}}{\Gamma(n+1)}\,\,\,n=0,1,2,\cdots$ be in set A.
Then its Shehu transform is given by
\begin{equation}
{\mathbb
S}\left[\frac{t^{n}}{\Gamma(n+1)}\right]=\left(\frac{u}{s}\right)^{n+1}.
\end{equation}
\end{pro}
The proof of property 6 follows immediately from the previous
property 5. \qquad\qquad $\Box$
\begin{pro}
Let the function $v(t)=\exp(\alpha t)$ be in A. Then its Shehu
transform is given by
\begin{equation}
{\mathbb S}\left[\exp(\alpha t)\right]=\frac{u}{s-\alpha u}.
\end{equation}
\textbf{Proof.} Using Equ.(\ref{eq2.1}), we get
\begin{eqnarray*}
{\mathbb S}\left[\exp(\alpha
t)\right]&=&\int_{0}^{\infty}\exp\left(-\frac{(s-\alpha
u)t}{u}\right)dt\\&&=-\frac{u}{s-\alpha
u}\lim_{\gamma\rightarrow\infty}\left[\exp\left(-\frac{(s-\alpha
u)t}{u}\right)\right]^\gamma_0=\frac{u}{s-\alpha u}\nonumber.
\end{eqnarray*}
This ends the proof.\qquad\qquad $\Box$
\end{pro}
\begin{pro}
Let the function $v(t)=t\exp(\alpha t)$ be in set A. Then its
Shehu transform is given by
\begin{equation}
{\mathbb S}\left[t\exp(\alpha t)\right]=\frac{u^2}{(s-\alpha
u)^2}.
\end{equation}
\vskip 1mm \textbf{Proof.} Using the Definition 1 of the Shehu
transform and integration by parts, we get
\begin{eqnarray*}
&&{\mathbb S}\left[t\exp(\alpha
t)\right]=\int_{0}^{\infty}t\exp\left(-\frac{(s-\alpha
u)t}{u}\right)dt\\&&=-\frac{u}{s-\alpha
u}\lim_{\gamma\rightarrow\infty}\left[t\exp\left(-\frac{(s-\alpha
u)t}{u}\right)\right]^\gamma_0\\&&+\frac{u}{s-\alpha
u}\int_{0}^{\infty}\exp\left(-\frac{(s-\alpha
u)t}{u}\right)dt\\&&=-\frac{u^2}{(s-\alpha
u)^2}\lim_{\gamma\rightarrow\infty}\left[\exp\left(-\frac{(s-\alpha
u)t}{u}\right)\right]^\gamma_0=\frac{u^2}{(s-\alpha u)^2}.
\end{eqnarray*}
The proof is complete.\qquad\qquad $\Box$
\end{pro}
\begin{pro}
Let the function $v(t)=\frac{t^{n}\exp(\alpha
t)}{n!}\,\,\,n=0,1,2,...$ be in set A. Then its Shehu transform is
given by
\begin{equation}
{\mathbb S}\left[\frac{t^{n}\exp(\alpha
t)}{n!}\right]=\frac{u^{n+1}}{(s-\alpha u)^{n+1}}.
\end{equation}
\vskip 1mm \textbf{Proof.} Using the Definition 1 of the Shehu
transform and integration by parts, we deduce
\begin{eqnarray*}
&&{\mathbb S}\left[t^{n}\exp(\alpha
t)\right]=\int_{0}^{\infty}t^{n}\exp\left(-\frac{(s-\alpha
u)t}{u}\right)dt\\&&=\frac{un}{(s-\alpha
u)}\int_{0}^{\infty}t^{n-1}\exp\left(-\frac{(s-\alpha
u)t}{u}\right)dt\\&=&\frac{u^2n(n-1)}{(s-\alpha
u)^2}\int_{0}^{\infty}t^{n-2}\exp\left(-\frac{(s-\alpha
u)t}{u}\right)dt=\cdots=\frac{n!}{(s-\alpha u)^{n+1}}.
\end{eqnarray*}
Thus the proof is complete.\qquad\qquad $\Box$
\end{pro}
\begin{pro}
Let the function $v(t)=\frac{t^{n}}{\Gamma(n+1)}\exp(\alpha
t)\,\,\,n=0,1,2,...$ be in set A. Then its Shehu transform is
given by
\begin{equation}
{\mathbb S}\left[\frac{t^{n}\exp(\alpha
t)}{\Gamma(n+1)}\right]=\frac{u^{n+1}}{(s-\alpha u)^{n+1}}.
\end{equation}
The proof of Property 10 follows as a direct consequence of
Property 9. \qquad\qquad $\Box$
\end{pro}
\begin{pro}
Let the function $v(t)=sin(\alpha t)$ be in set A. Then its Shehu
transform is given by
\begin{equation}
{\mathbb S}\left[\sin(\alpha t)\right]=\frac{\alpha
u^{2}}{s^2+\alpha^2u^2}.
\end{equation}
\textbf{Proof.} Using the Definition 1 of the Shehu transform and
integration by parts, we get
\begin{eqnarray*}
&&{\mathbb S}\left[\sin(\alpha
t)\right]=\int_{0}^{\infty}\exp\left(-\frac{st}{u}\right)\sin(\alpha
t)dt\\&&=-\frac{
u}{s}\lim_{\gamma\rightarrow\infty}\left[\exp\left(-\frac{st}{u}\right)\sin(\alpha
t)\right]^\gamma_0+\frac{
u\alpha}{s}\int_{0}^{\infty}\exp\left(-\frac{st}{u}\right)\cos(\alpha
t)dt\\&=&-\frac{\alpha
u^2}{s^2}\lim_{\gamma\rightarrow\infty}\left[\exp\left(-\frac{st}{u}\right)\cos(\alpha
t)\right]^\gamma_0-\frac{\alpha^2u^2}{s^2}\int_{0}^{\infty}\exp\left(-\frac{st}{u}\right)\sin(\alpha
t)dt\\&=&\frac{\alpha
u^2}{s^2}-\frac{\alpha^2u^2}{s^2}\int_{0}^{\infty}\exp\left(-\frac{st}{u}\right)\sin(\alpha
t)dt.
\end{eqnarray*}
Simplifying the required integrals complete the proof of Property
11. \qquad\qquad $\Box$
\end{pro}
\begin{pro}
Let the function $v(t)=\cos(\alpha t)$ be in set A. Then its Shehu
transform is given by
\begin{equation}
{\mathbb S}\left[\cos(\alpha t)\right]=\frac{us}{s^2+\alpha^2u^2}.
\end{equation}
\textbf{Proof.} Using the Definition 1 of the Shehu transform and
integration by parts, we deduce
\begin{eqnarray*}
&&{\mathbb S}\left[\cos(\alpha
t)\right]=\int_{0}^{\infty}\exp\left(-\frac{st}{u}\right)\cos(\alpha
t)dt\\&&=-\frac{u}{s}\lim_{\gamma\rightarrow\infty}\left[\exp\left(-\frac{st}{u}\right)\cos(\alpha
t)\right]^\gamma_0-\frac{\alpha
u}{s}\int_{0}^{\infty}\exp\left(-\frac{st}{u}\right)\sin(\alpha
t)dt\\&=&\frac{u}{s}-\frac{\alpha
u^2}{s^2}\lim_{\gamma\rightarrow\infty}\left[\exp\left(-\frac{st}{u}\right)\sin(\alpha
t)\right]^\gamma_0-\frac{\alpha^2u^2}{s^2}\int_{0}^{\infty}\exp\left(-\frac{st}{u}\right)\cos(\alpha
t)dt\\&=&\frac{u}{s}-\frac{\alpha^2u^2}{s^2}\int_{0}^{\infty}\exp\left(-\frac{st}{u}\right)\cos(\alpha
t)dt.
\end{eqnarray*}
Simplifying the required integrals complete the proof of Property
12. \qquad\qquad $\Box$
\end{pro}
\begin{pro}
Let the function $v(t)=\frac{\sinh(\alpha t)}{\alpha}$ be in set
A. Then its Shehu transform is given by
\begin{equation}
{\mathbb S}\left[\frac{\sinh(\alpha t)}{\alpha}\right]=\frac{
u^{2}}{s^2-\alpha^2u^2}.
\end{equation}
\textbf{Proof.} From the Definition 1 of the Shehu transform and
integration by parts, we get
\begin{eqnarray*}
&&{\mathbb S}\left[\sinh(\alpha
t)\right]=\int_{0}^{\infty}\exp\left(-\frac{st}{u}\right)\sinh(\alpha
t)dt\\&&=-\frac{
u}{s}\lim_{\gamma\rightarrow\infty}\left[\exp\left(-\frac{st}{u}\right)\sinh(\alpha
t)\right]^\gamma_0+\frac{
u\alpha}{s}\int_{0}^{\infty}\exp\left(-\frac{st}{u}\right)\cosh(\alpha
t)dt\\&=&-\frac{\alpha
u^2}{s^2}\lim_{\gamma\rightarrow\infty}\left[\exp\left(-\frac{st}{u}\right)\cos(\alpha
t)\right]^\gamma_0+\frac{\alpha^2u^2}{s^2}\int_{0}^{\infty}\exp\left(-\frac{st}{u}\right)\sinh(\alpha
t)dt\\&=&\frac{\alpha
u^2}{s^2}+\frac{\alpha^2u^2}{s^2}\int_{0}^{\infty}\exp\left(-\frac{st}{u}\right)\sinh(\alpha
t)dt.
\end{eqnarray*}
Simplifying the required integrals complete the proof of Property
13. \qquad\qquad $\Box$
\end{pro}
\begin{pro}
Let the function $v(t)=\cosh(\alpha t)$ be in set A. Then its
Shehu transform is given by
\begin{equation}
{\mathbb S}\left[\cosh(\alpha
t)\right]=\frac{us}{s^2-\alpha^2u^2}.
\end{equation}
\textbf{Proof.} Applying the Definition 1 of the Shehu transform
and integration by parts, we get
\begin{eqnarray*}
&&{\mathbb S}\left[\cosh(\alpha
t)\right]=\int_{0}^{\infty}\exp\left(-\frac{st}{u}\right)\cosh(\alpha
t)dt\\&&=-\frac{u}{s}\lim_{\gamma\rightarrow\infty}\left[\exp\left(-\frac{st}{u}\right)\cos(\alpha
t)\right]^\gamma_0+\frac{\alpha
u}{s}\int_{0}^{\infty}\exp\left(-\frac{st}{u}\right)\sinh(\alpha
t)dt\\&=&\frac{u}{s}-\frac{\alpha
u^2}{s}\lim_{\gamma\rightarrow\infty}\left[\exp\left(-\frac{st}{u}\right)\sinh(\alpha
t)\right]^\gamma_0+\frac{\alpha^2u^2}{s^2}\int_{0}^{\infty}\exp\left(-\frac{st}{u}\right)\cos(\alpha
t)dt\\&=&\frac{u}{s}+\frac{\alpha^2u^2}{s^2}\int_{0}^{\infty}\exp\left(-\frac{st}{u}\right)\cosh(\alpha
t)dt.
\end{eqnarray*}
Collecting the required integrals complete the proof of Property
14. \qquad\qquad $\Box$
\end{pro}
\begin{pro}
Let the function $\frac{\exp\left(\beta t\right)\sin(\alpha
t)}{\alpha}$ be in set A. Then its Shehu transform is given by
\begin{equation}
{\mathbb S}\left[\frac{\exp\left(\beta t\right)\sin(\alpha
t)}{\alpha}\right]=\frac{u^2}{(s-\beta u)^2+\alpha^2u^2}.
\end{equation}
\textbf{Proof.} Using the Definition 1 of the Shehu transform and
integration by parts, we deduce
\begin{eqnarray*}
{\mathbb S}\left[\exp\left(\beta t\right)\sin(\alpha
t)\right]&=&\int_{0}^{\infty}\exp\left(-\frac{(s-\beta
u)}{u}t\right)\sin(\alpha t)dt\\&&=\frac{-u}{(s-\beta
u)}\lim_{\gamma\rightarrow\infty}\left[\exp\left(-\frac{(s-\beta
u)}{u}t\right)\sin(\alpha
t)dt\right]^\gamma_0\\
&&+\frac{u\alpha}{s-\beta
u}\int_{0}^{\infty}\exp\left(-\frac{(s-\beta
u)}{u}t\right)\cos(\alpha t)dt\\&&=-\frac{u^2\alpha}{(s-\beta
u)^2}\lim_{\gamma\rightarrow\infty}\left[\exp\left(-\frac{(s-\beta
u)}{u}t\right)\cos(\alpha t)dt\right]^\gamma_0\\&&-\frac{\alpha
u^2\alpha^2}{(s-\beta
u)^2}\int_{0}^{\infty}\exp\left(-\frac{(s-\beta
u)}{u}t\right)\sin(\alpha t)dt\\&=&\frac{u^2\alpha}{(s-\beta
u)^2}-\frac{ u^2\alpha^2}{(s-\beta
u)^2}\int_{0}^{\infty}\exp\left(-\frac{(s-\beta
u)}{u}t\right)\sin(\alpha t)dt.
\end{eqnarray*}
Simplifying the required integrals complete the proof of property
15. This ends the proof.\qquad\qquad $\Box$
\end{pro}
\begin{pro}
Let the function $\exp\left(\beta t\right)\cos(\alpha t)$ be in
set A. Then its Shehu transform is given by
\begin{equation}
{\mathbb S}\left[\exp\left(\beta t\right)\cos(\alpha
t)\right]=\frac{u(s-\alpha u)}{(s-\beta u)^2+\alpha^2u^2}.
\end{equation}
\textbf{Proof.} Applying the Definition 1 of the Shehu transform
and integration by parts, we get
\begin{eqnarray*}
{\mathbb S}\left[\exp\left(\beta t\right)\cos(\alpha
t)\right]&=&\int_{0}^{\infty}\exp\left(-\frac{(s-\beta
u)}{u}t\right)\cos(\alpha t)dt\\&&=-\frac{u}{s-\beta
u}\lim_{\gamma\rightarrow\infty}\left[\exp\left(-\frac{(s-\beta
u)}{u}t\right)\cos(\alpha t)\right]^\gamma_0\\&&+\frac{\alpha
u}{s-\beta u}\int_{0}^{\infty}\exp\left(-\frac{(s-\beta
u)}{u}t\right)\sin(\alpha t)dt\\ &=&\frac{u}{s-\beta
u}+\frac{\alpha u}{s-\beta
u}\int_{0}^{\infty}\exp\left(-\frac{(s-\beta
u)}{u}t\right)\sin(\alpha t)dt\\&=&\frac{u}{s-\beta
u}+\frac{\alpha u^2}{(s-\alpha
u)^2}\lim_{\gamma\rightarrow\infty}\left[\exp\left(-\frac{(s-\beta
u)}{u}t\right)\sin(\alpha t)\right]^\gamma_0\\&&-\frac{\alpha^2
u^2}{(s-\beta u)^2}\int_{0}^{\infty}\exp\left(-\frac{(s-\beta
u)}{u}t\right)\cos(\alpha t)dt\\&=&\frac{u}{s-\beta
u}-\frac{\alpha^2 u^2}{(s-\beta
u)^2}\int_{0}^{\infty}\exp\left(-\frac{(s-\beta
u)}{u}t\right)\cos(\alpha t)dt.
\end{eqnarray*}
Simplifying the required integrals complete the proof of property
16.\qquad\qquad $\Box$
\end{pro}
\begin{pro}
Let the function $\frac{\exp\left(\beta t\right)-\exp\left(\alpha
t\right)}{\beta-\alpha}$ be in set A. Then its Shehu transform is
given by
\begin{equation}
{\mathbb S}\left[\frac{\exp\left(\alpha
t\right)}{\beta-\alpha}\right]=\frac{u^2}{(s-\alpha u)(s-\beta
u)}.
\end{equation}
\textbf{Proof.} Using the definition of Shehu transform, we get
\begin{eqnarray*}
&&{\mathbb S}\left[\frac{\exp\left(\alpha
t\right)}{\beta-\alpha}\right]=\frac{u}{\beta-\alpha}\int_{0}^{\infty}\exp\left(\frac{-st}{u}\right)\left(\exp(\beta
t)-\exp\left(\alpha
t\right)\right)dt\\&&=\frac{1}{\beta-\alpha}\int_{0}^{\infty}e^{\frac{(\beta
u-s)}{u}t}dt-\frac{1}{\beta-\alpha}\int_{0}^{\infty}\exp\left(\frac{(\alpha
u-s)t}{u}\right)dt\\&&=\frac{u}{(\beta-\alpha)(\beta
u-s)}\lim_{\gamma\rightarrow\infty}\left[\exp\left(-\frac{(s-\beta
u)t}{u}\right)\right]^\gamma_0\\&&-\frac{u}{(\beta-\alpha)(\alpha
u-s)}\lim_{\gamma\rightarrow\infty}\left[\exp\left(-\frac{(s-\beta
u)t}{u}\right)\right]^\gamma_0\\&=&-\frac{u}{(\beta-\alpha)(\beta
u-s)}+\frac{u}{(\beta-\alpha)(\alpha u-s)}\\&=&\frac{-u(\alpha
u-s)+u(\beta u-s)}{(\beta-\alpha)(\alpha u-s)(\beta
u-s)}=\frac{u^2}{(s-\alpha u)(s-\beta u)}.
\end{eqnarray*}
The proof is complete.\qquad\qquad $\Box$
\end{pro}
\begin{pro}
Let the function $\frac{\beta \exp\left(\beta t\right)-\alpha
\exp(\alpha t)}{\beta-\alpha}$ be in set A. Then its Shehu
transform is given by:
\begin{equation}
{\mathbb S}\left[\frac{\beta \exp\left(\beta t\right)-\alpha
\exp\left(\alpha
t\right)}{\beta-\alpha}\right]=\frac{us}{(s-\alpha u)(s-\beta u)}.
\end{equation}
\textbf{Proof:}\\
Using the definition of Shehu transform, we get
\begin{eqnarray*}
&&{\mathbb S}\left[\frac{\beta \exp\left(\beta t\right)-\alpha
\exp\left(\alpha
t\right)}{\beta-\alpha}\right]=\frac{1}{\beta-\alpha}\int_{0}^{\infty}\exp\left(\frac{-st}{u}\right)\left(\beta
\exp\left(\beta t\right)-\alpha \exp\left(\alpha
t\right)\right)dt\\&&=\frac{\beta}{\beta-\alpha}\int_{0}^{\infty}\exp\left(\frac{(\beta
u-s)}{u}t\right)dt-\frac{\alpha}{\beta-\alpha}\int_{0}^{\infty}\exp\left(\frac{(\alpha
u-s)t}{u}\right)dt\\&&=\frac{u\beta}{(\beta-\alpha)(\beta
u-s)}\lim_{\gamma\rightarrow\infty}\left[\exp\left(-\frac{(s-\beta
u)t}{u}\right)\right]^\gamma_0\\&&-\frac{u\alpha}{(\beta-\alpha)(\alpha
u-s)}\lim_{\gamma\rightarrow\infty}\left[\exp\left(-\frac{(s-\alpha
u)t}{u}\right)\right]^\gamma_0\\&=&-\frac{u\beta}{(\beta-\alpha)(\beta
u-s)}+\frac{u\alpha}{(\beta-\alpha)(\alpha
u-s)}\\&=&\frac{-u\beta(\alpha u-s)+u\alpha(\beta
u-s)}{(\beta-\alpha)(\alpha u-s)(\beta u-s)}=\frac{us}{(s-\alpha
u)(s-\beta u)}.
\end{eqnarray*}
This ends the proof.\qquad\qquad $\Box$
\end{pro}
More properties of the Shehu transform and their converges to the
natural transform, the Sumudu transform, and the Laplace transform
are presented in table 1. The comprehensive summary of Shehu
transform properties are presented in table 2.


\section{Applications}
In this section, the applications of the proposed transform are
presented. The simplicity, efficiency and high accuracy of the
Shehu transform are clearly illustrated.
\begin{exa}
Consider the following first order ordinary differential equation
\begin{equation}\label{eq4.1}
\frac{dv(t)}{dt}+v(t)=0,
\end{equation}
subject to the initial condition
\begin{equation}
v(0)=1.
\end{equation}
Applying the Shehu transform on both sides of Equ. (\ref{eq4.1}),
we get
\begin{equation}
\frac{s}{u}V(s,u)-v(0)+V(s,u)=0.
\end{equation}
\end{exa}
Substituting the given initial condition and simplifying, we
deduce
\begin{equation}\label{eq4.4}
V(s,u)=\frac{u}{s+u}
\end{equation}
Taking the inverse Shehu transform of Equ. (\ref{eq4.4}), yields
\begin{equation}
v(t)=\exp(-t).
\end{equation}
\begin{exa}
Consider the following second order ordinary differential equation
\begin{equation}\label{eq4.6}
\frac{d^2v(t)}{dt^2}+\frac{dv(t)}{dt}=1
\end{equation}
subject to the initial conditions
\begin{equation}
v(0)=0,\,\,\,\,\frac{dv(0)}{dt}=0.
\end{equation}
Applying the Shehu transform on both sides of Equ. (\ref{eq4.6}),
we obtain
\begin{equation}
\frac{s^2}{u^2}V(s,u)-\frac{s}{u}v(0)-v'(0)+\frac{s}{u}V(s,u)-v(0)=\frac{u}{s}.
\end{equation}
\end{exa}
Substituting the given initial conditions and simplifying, we
deduce
\begin{equation}\label{eq4.9}
V(s,u)=-\frac{u}{s}+\frac{u^2}{s^2}+\frac{u}{s+u}.
\end{equation}
Taking the inverse Shehu transform of Equ. (\ref{eq4.9}), we get
\begin{equation}
v(t)=-1+t+\exp(-t).
\end{equation}

\begin{exa}
Consider the following second nonhomogeneous order ordinary
differential equation
\begin{equation}\label{eq4.11}
\frac{d^2v(t)}{dt^2}-3\frac{dv(t)}{dt}+2v(t)=\exp(3t)
\end{equation}
subject to the initial conditions
\begin{equation}
v(0)=1,\,\,\,\,\frac{dv(0)}{dt}=0.
\end{equation}
Applying the Shehu transform on both sides of Equ. (\ref{eq4.11}),
yields
\begin{equation}
\frac{s^2}{u^2}V(s,u)-\frac{s}{u}v(0)-v'(0)-3\left(\frac{s}{u}V(s,u)-v(0)\right)+2V(s,u)=\frac{u}{s-3u}.
\end{equation}
\end{exa}
Substituting the given initial conditions and simplifying, we
obtain
\begin{equation}\label{eq4.14}
V(s,u)=\frac{5}{2}\frac{u}{(s-u)}-2\frac{u}{s-2u}+\frac{1}{2}\frac{u}{(s-3u)}.
\end{equation}
Taking the inverse Shehu transform of Equ. (\ref{eq4.14}), we get
\begin{equation}
v(t)=\frac{5}{2}\exp(t)-2\exp(2t)+\frac{1}{2}\exp(3t).
\end{equation}

\begin{exa}
Consider the following ordinary differential equation
\begin{equation}\label{eq4.16}
\frac{d^2v(t)}{dt^2}+2\frac{dv(t)}{dt}+5v(t)=\exp(-t)\sin (t)
\end{equation}
subject to the initial conditions
\begin{equation}
v(0)=0,\,\,\,\,\frac{dv(0)}{dt}=1.
\end{equation}
Applying the Shehu transform on both sides of Equ. (\ref{eq4.16}),
we get
\begin{equation}
\frac{s^2}{u^2}V(s,u)-\frac{s}{u}v(0)-v'(0)+2\left(\frac{s}{u}V(s,u)-\frac{s}{u}v(0)\right)+5V(s,u)=\frac{u^2}{(s+u)^2+u^2}
\end{equation}
\end{exa}
Substituting the given initial conditions and simplifying, we get
\begin{equation}\label{eq4.19}
V(s,u)=\frac{1}{3}\frac{u^2}{((s+u)^2+u^2)}+\frac{2}{3}\frac{u^2}{((s+u)^2+(2u)^2)}
\end{equation}
Taking the inverse Shehu transform of Equ. (\ref{eq4.19}), we get
\begin{equation}
v(t)=\frac{1}{3}\exp(-t)\sin (t)+\frac{2}{3}\exp(-t)\sin (2t).
\end{equation}
\begin{exa}
Consider the following homogeneous partial differential equation
\begin{equation}\label{eq4.21}
\frac{\partial v(x,t)}{\partial t}=\frac{\partial^2
v(x,t)}{\partial x^2}
\end{equation}
subject to the boundary and initial conditions
\begin{equation}
v(0,t)=0,\,\,\,v(1,t)=0,\,\,\,\,v(x,0)=3sin(2\pi x).
\end{equation}
Applying the Shehu transform on both sides of Equ. (\ref{eq4.21}),
we get
\begin{equation}
\frac{s}{u}V(x,s,u)-v(x,0)=\frac{d^2V(x,s,u)}{dx^2}.
\end{equation}
\end{exa}
Substituting the given initial condition and simplifying, we get
\begin{equation}\label{eq4.24}
\frac{d^2V(x,s,u)}{dx^2}-\frac{s}{u}V(x,s,u)=-3sin(2\pi x).
\end{equation}
The general solution of Equ. (\ref{eq4.24}) can be written as
\begin{equation}\label{eq4.25}
V(x,s,u)=V_h(x,s,u)+V_p(x,s,u),
\end{equation}
where $V_h(x,s,u)$ is the solution of the homogeneous part which
is given by
\begin{equation}\label{eq4.26}
V_h(x,s,u)=\alpha_1\exp\left(\sqrt{\frac{s}{u}}x\right)+\alpha_2\exp\left(-\sqrt{\frac{s}{u}}x\right),
\end{equation}
and $V_p(x,s,u)$ is the solution of the nonhomogeneous part which
is given by
\begin{equation}\label{eq4.27}
V_p(x,s,u)=\beta_1sin(2\pi x)+\beta_2cos(2\pi x).
\end{equation}
Applying the boundary conditions on Equ. (\ref{eq4.26}), we get
\begin{equation}
\alpha_1+\alpha_2=0\,\,\,\,and\,\,\,\,
\alpha_1\exp\left(\sqrt{\frac{s}{u}}\right)+\alpha_2\exp\left(-\sqrt{\frac{s}{u}}\right)=0\Rightarrow
V_h(x,s,u)=0,\nonumber
\end{equation}
 since $\alpha_1=\alpha_2=0$.\\
Using the method of undetermined coefficients on the
nonhomogeneous part, we get
\begin{equation}
V_p(x,s,u)=\frac{3u}{s+4\pi^2 u}\sin(2\pi x),
\end{equation}
Since, $\beta_1=\frac{3u}{s+4\pi^2
u},$\,\,\,\,and\,\,\,\,$\beta_2=0$.\\
Then Equ. (\ref{eq4.25}) will become
\begin{equation}\label{eq4.29}
V(x,s,u)=\frac{3u}{s+4\pi^2 u}\sin(2\pi x),
\end{equation}
 Taking the inverse Shehu transform of Equ. (\ref{eq4.29}),
we get
\begin{equation}
v(x,t)=3\exp(-4\pi^2t)\sin(2\pi x).
\end{equation}
\begin{figure}[!h]
    \centering
    \begin{minipage}[t]{2cm}
        \centering
        \includegraphics[scale=0.4]{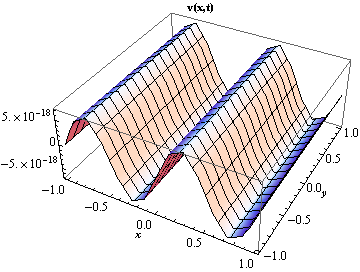}
    \end{minipage}
    \hspace{3cm}
    \begin{minipage}[t]{2cm}
        \centering
        \includegraphics[scale=0.4]{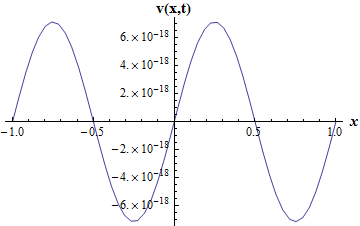}
    \end{minipage}
    \end{figure}\\
    \textbf{Figure 1.} 3D and 2D surfaces of the analytical solution of Equ. (\ref{eq4.21}) in the ranges $-1<x<1$, and $-1<t<1$, when $t=1$.
\begin{exa}
Consider the following nonhomogeneous partial differential
equation
\begin{equation}\label{eq4.31}
\frac{\partial^2 v(x,t)}{\partial t^2}=\beta^2\frac{\partial^2
v(x,t)}{\partial x^2}+\sin(\pi x)
\end{equation}
subject to the boundary and initial conditions
\begin{equation}
v(0,t)=0,\,\,\,v(1,t)=0,\,\,\,\,v(x,0)=0,\,\,\,\frac{\partial
v(x,0)}{\partial t}=0,\,\,\,\beta^2=1.
\end{equation}
Applying the Shehu transform on both sides of Equ. (\ref{eq4.31}),
we get
\begin{equation}
\frac{s^2}{u^2}V(x,s,u)-\frac{s}{u}v(x,0)-v'(x,0)=\frac{d^2V(x,s,u)}{dx^2}+\frac{u}{s}\sin(\pi
x).
\end{equation}
\end{exa}
Substituting the given initial condition and simplifying, we get
\begin{equation}\label{eq4.34}
\frac{d^2V(x,s,u)}{dx^2}-\frac{s^2}{u^2}V(x,s,u)=-\frac{u}{s}\sin(\pi
x).
\end{equation}
The general solution of Equ. (\ref{eq4.34}) can be written as
\begin{equation}\label{eq4.35}
V(x,s,u)=V_h(x,s,u)+V_p(x,s,u),
\end{equation}
where $V_h(x,s,u)$ is the solution of the homogeneous part which
is given by
\begin{equation}\label{eq4.36}
V_h(x,s,u)=\lambda_1\exp\left(\frac{s}{u}x\right)+\lambda_2\exp\left(-\frac{s}{u}x\right),
\end{equation}
and $V_p(x,s,u)$ is the solution of the nonhomogeneous part which
is given by
\begin{equation}\label{eq4.37}
V_p(x,s,u)=\eta_1\sin(\pi x)+\eta_2\cos(\pi x).
\end{equation}
Applying the boundary conditions on Equ. (\ref{eq4.36}), we deduce
\begin{equation}
\lambda_1+\lambda_2=0\,\,\,\,and\,\,\,\,
\lambda_1\exp\left(\frac{s}{u}\right)+\lambda_2\exp\left(-\frac{s}{u}\right)=0\Rightarrow
V_h(x,s,u)=0,\nonumber
\end{equation}
 since $\lambda_1=\lambda_2=0$.\\
Using the method of undetermined coefficients on the
nonhomogeneous part, we get
\begin{equation}
V_p(x,s,u)=\frac{1}{\pi^2}\left(\frac{u}{s}-\frac{us}{s^2+u^2\pi^2}\right)\sin(\pi
x),
\end{equation}
since,
\begin{equation}
\eta_1=\frac{u^3}{s(s^2+u^2
\pi^2)}=\frac{1}{\pi^2}\left(\frac{u}{s}-\frac{us}{s^2+u^2\pi^2}\right),\,\,\,and\,\,\,\,\eta_2=0.\nonumber
\end{equation}
Then Equ. (\ref{eq4.35}) will becomes
\begin{equation}\label{eq4.39}
V(x,s,u)=\frac{1}{\pi^2}\left(\frac{u}{s}-\frac{us}{s^2+u^2\pi^2}\right)\sin(\pi
x).
\end{equation}
 Taking the inverse Shehu transform of Equ. (\ref{eq4.39}),
we get
\begin{equation}
v(x,t)=\frac{1}{\pi^2}\left(1-\cos(\pi t)\right)\sin(\pi x).
\end{equation}
\begin{figure}[!h]
    \centering
    \begin{minipage}[t]{2cm}
        \centering
        \includegraphics[scale=0.4]{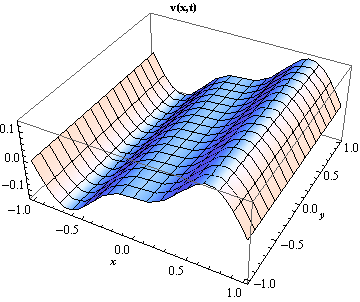}
    \end{minipage}
    \hspace{3cm}
    \begin{minipage}[t]{2cm}
        \centering
        \includegraphics[scale=0.4]{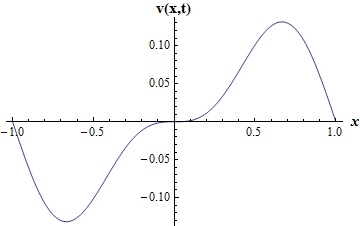}
    \end{minipage}
    \end{figure}\\
    \textbf{Figure 2.} 3D and 2D surfaces of the analytical solution of Equ. (\ref{eq4.31}) in the ranges $-1<x<1$, and $-1<t<1$.
\begin{center}\section{Conclusion}\end{center}
We introduced an efficient Laplace-type integral transform called
the Shehu transform for solving both ordinary and partial
differential equations. We presented its existence and inverse
transform. We presented some useful properties and their
applications for solving ordinary and partial differential
equations. We provide a comprehensive list of the Laplace
transform, Sumudu transform, and the natural transform in table 1
to show their mutual relationship with the Shehu transform.
Finally, based on the mathematical formulations, simplicity and
the findings of the proposed integral transform, we conclude that
it is highly efficient because of the following advantages:
\begin{itemize}
   \item It is a generalization of the Laplace and the Sumudu integral transforms.
   \item Its visualization is easier than the Sumudu transform, the natural transform, and the Elzaki transform.
   \item The Laplace-type integral transform become Laplace transform when the variable $u=1$ and the Yang integral transform when the variable $s=1$.
   \item It can easily be applied directly to some class of ordinary and the partial differential equations as demonstrated in the application section.
   \item For advanced research in physical science and engineering, the proposed integral transform can be considered a stepping-stone to the Sumudu transform, the natural transform, the Elzaki transform, and the Laplace transform.
\end{itemize}

\section{Acknowledgements}
The authors are highly grateful to the editor's and the anonymous
referees' for their useful comments and suggestions in this paper.
This research is partially supported by the National Natural
Science Foundations of China under Grants No. (11571206). The
first author also acknowledges the financial support of China
Scholarship Council (CSC) in Shandong University with grand (CSC
No: 2017GXZ025381).

\small \begin{center}
\end{center}
\section*{Appendix}
\begin{center}
\textbf{Table 1:} Here we present a comprehensive list of the
Shehu transform of some special functions and their relationship
with the natural transform ${\mathbb N}\left[v(t)\right]$, the
Sumudu transform $\textbf{S}\left[v(t)\right]$, and the Laplace
transform. {\setlength{\extrarowheight}{10pt}
\begin{longtable}{|c|c|c|c|c|c|}
  \hline
  S.NO. & $v(t)$ & ${\mathbb S}\left[v(t)\right]$ & ${\mathbb N}\left[v(t)\right]$ & $\textbf{S}\left[v(t)\right]$ & $\pounds\left[v(t)\right]$ \\
  \hline
  1 & 1 &$\frac{u}{s}$& $\frac{1}{s}$& 1 & $\frac{1}{s}$ \\
  2 & $t$& $\frac{u^2}{s^2}$& $\frac{u}{s^2}$& $u$& $\frac{1}{s^2}$\\
  3 & $\exp(\alpha (t))$ & $\frac{u}{s-\alpha u}$ & $\frac{1}{s-\alpha u}$ & $\frac{1}{1-\alpha u}$ & $\frac{1}{s-\alpha}$\\
  4 & $\frac{\sin(\alpha t)}{\alpha}$ & $\frac{u^2}{s^2+\alpha^2u^2}$ & $\frac{u}{s^2+\alpha^2u^2}$ & $\frac{u}{1+\alpha^2u^2}$ & $\frac{1}{s^2+\alpha^2}$ \\
  5 & $\cos(\alpha t)$ & $\frac{us}{s^2+\alpha^2u^2}$ & $\frac{s}{s^2+\alpha^2u^2}$ & $\frac{1}{1+\alpha^2u^2}$ & $\frac{s}{s^2+\alpha^2}$ \\
  6 & $\cosh\alpha t$ & $\frac{us}{s^2-u^2}$ & $\frac{s}{s^2-u^2}$ & $\frac{1}{1-u^2}$ & $\frac{s}{s^2-1}$ \\
  7 & $\frac{t^{n}}{n!}\,\,n=0,1,2...$ & $\left(\frac{u}{s}\right)^{n+1}$ & $\frac{u^{n}}{s^{n+1}}$ & $u^{n}$ & $\frac{1}{s^{n+1}}$ \\
  8 & $\frac{t^{n}}{\Gamma(n+1)}\,\,n=0,1,2,...$ & $\left(\frac{u}{s}\right)^{n+1}$ & $\frac{u^{n}}{s^{n+1}}$ & $u^{n}$ & $\frac{1}{s^{n+1}}$ \\
  9 & $\cos(t)$ & $\frac{us}{s^2+u^2}$ & $\frac{s}{s^2+u^2}$ & $\frac{1}{1+u^2}$ & $\frac{1}{s^2+1}$ \\
  10 & $\sin(t)$ & $\frac{u^2}{s^2+u^2}$ & $\frac{u}{s^2+u^2}$ & $\frac{u}{1+u^2}$ & $\frac{1}{s^2+1}$ \\
  11 & $\frac{sinh(\alpha t)}{\alpha}$ & $\frac{u^2}{s^2-\alpha^2u^2}$ & $\frac{\alpha u}{s^2-\alpha^2u^2}$ & $\frac{\alpha u^2}{1-\alpha^2u^2}$ & $\frac{\alpha}{s^2-\alpha^2}$ \\
  12 & $\cosh(\alpha t)$ & $\frac{us}{s^2-\alpha^2u^2}$ & $\frac{s}{s^2-\alpha^2u^2}$ & $\frac{1}{1-\alpha^2u^2}$ & $\frac{s}{s^2-\alpha^2}$ \\
  13 & $\exp(\beta t)cosh(\alpha t)$ & $\frac{u(s-\beta u)}{(s-\beta u)^2-\alpha^2u^2}$ & $\frac{s-\beta u}{(s-\beta u)^2-\alpha^2u^2}$ & $\frac{1-\beta u}{(s-\beta u)^2-\alpha^2u^2}$ & $\frac{s-\beta }{(s-\beta u)^2-\alpha^2}$ \\
  14 & $\frac{\exp(\beta t)\sinh(\alpha t)}{\alpha}$ & $\frac{u^2}{(s-\beta u)^2-\alpha^2u^2}$ & $\frac{u}{(s-\beta u)^2-\alpha^2u^2}$ & $\frac{u}{(1-\beta u)^2-\alpha^2u^2}$ & $\frac{1}{(s-\beta )^2-\alpha^2}$ \\
  15 & $\frac{t\sin(\alpha t)}{2\alpha}$ & $\frac{u^3s}{(s^2+\alpha^2u^2)^2}$ & $\frac{u^2s}{(s^2+\alpha^2u^2)^2}$ & $\frac{u^3}{(1+\alpha^2u^2)^2}$ & $\frac{s}{(s^2+\alpha^2)^2}$ \\
  16 & $t\cos(\alpha t)$ & $\frac{u^2(s^2-\alpha^2u^2)^2}{(s^2+\alpha^2u^2)^2}$ & $\frac{u(s^2-\alpha^2u^2)^2}{(s^2+\alpha^2u^2)^2}$ & $\frac{u(1-\alpha^2u^2)^2}{(1+\alpha^2u^2)^2}$ & $\frac{(s^2-\alpha^2)^2}{(s^2+\alpha^2)^2}$ \\
17 & $\frac{\sin(\alpha t)+\alpha t\cos(\alpha t)}{2\alpha}$ & $\frac{u^2s^2}{(s^2+\alpha^2u^2)^2}$ & $\frac{us^2}{(s^2+\alpha^2u^2)^2}$ & $\frac{u}{(1+\alpha^2u^2)^2}$ & $\frac{s^2}{(s^2+\alpha^2)^2}$ \\
  18 & $\cos\alpha t-\frac{\alpha t\sin(\alpha t)}{2}$ & $\frac{us^3}{(s^2+\alpha^2u^2)^2}$ & $\frac{s^3}{(s^2+\alpha^2u^2)^2}$ & $\frac{1}{(1+\alpha^2u^2)^2}$ & $\frac{s^3}{(s^2+\alpha^2)^2}$ \\
  19 & $\frac{\sin(\alpha t)-\alpha t\cos(\alpha t)}{2\alpha^3}$ & $\frac{u^4}{(s^2+\alpha^2u^2)^2}$ & $\frac{u^3}{(s^2+\alpha^2u^2)^2}$ & $\frac{u^3}{(1+\alpha^2u^2)^2}$ & $\frac{1}{(s^2+\alpha^2)^2}$ \\
  20 & $t\sinh(\alpha t)+t\cosh(\alpha t)$ & $\frac{u^2}{(s-\alpha u)^2}$ & $\frac{u}{(s-\alpha u)^2}$ & $\frac{u^2}{(1-\alpha u)^2}$ & $\frac{1}{(s-\alpha)^2}$ \\
  21 & $\frac{t\sinh(\alpha t)}{2\alpha}$ & $\frac{u^3s}{(s^2-\alpha^2u^2)^2}$ & $\frac{u^2s}{(s^2-\alpha^2u^2)^2}$ & $\frac{u^2}{(1-\alpha^2u^2)^2}$ & $\frac{s}{(s^2-\alpha^2)^2}$ \\
  22 & $Si(\alpha t)$ (Sine integral) & $\frac{u}{s}tan^{-1}\left(\frac{\alpha u}{s}\right)$ & $\frac{1}{s}tan^{-1}\left(\frac{\alpha u}{s}\right)$ & $tan^{-1}\left(u\sqrt{\alpha^2}\right)$ & $\frac{1}{s}tan^{-1}\left(\frac{\alpha}{s}\right)$ \\
  23 & $Ci(\alpha t)$ (Cosine integral) & $-\frac{u}{2s}\log\left(\frac{s^2+\alpha^2}{\alpha^2}\right)$ & $-\frac{1}{2s}\log\left(\frac{s^2+\alpha^2u^2}{\alpha^2u^2}\right)$ & $-\frac{1}{2}\log\left(\frac{\alpha^2u^2+1}{\alpha^2u^2}\right)$ &$-\frac{1}{2s}\log\left(\frac{s^2+\alpha^2}{\alpha^2}\right)$ \\
  24 & $Ei(\alpha t)$ (Exp. integral)& $-\frac{u}{s}\log\left(\frac{\alpha u-s}{\alpha u}\right)$ & $-\frac{1}{s}\log\left(\frac{\alpha u-s}{\alpha u}\right)$ & $\log\left(\frac{\alpha u-1}{\alpha u}\right)$ &$-\frac{1}{s}\log\left(\frac{\alpha-s}{\alpha}\right)$\\
  25 & $\frac{(3-\alpha^2t^2)\sin(\alpha t)-3\alpha t\cos(\alpha t)}{8\alpha^5}$ & $\frac{u^6}{(s^2+\alpha^2u^2)^3}$ & $\frac{u^5}{(s^2+\alpha^2u^2)^3}$ & $\frac{u^5}{(1+\alpha^2u^2)^3}$ & $\frac{1}{(s^2+\alpha^2)^3}$ \\
  26 & $\frac{(3-\alpha^2t^2)\sin(\alpha t)+5\alpha t\cos(\alpha t)}{8\alpha}$& $\frac{u^2s^4}{(s^2+\alpha^2u^2)^3}$ & $\frac{us^4}{(s^2+\alpha^2u^2)^3}$ & $\frac{u}{(1+\alpha^2u^2)^3}$ & $\frac{s^4}{(s^2+\alpha^2)^3}$ \\
  27 & $\frac{(8-\alpha^2t^2)\cos(\alpha t)-7\alpha t\sin(\alpha t)}{8}$ & $\frac{us^5}{(s^2+\alpha^2u^2)^3}$ & $\frac{s^5}{(s^2+\alpha^2u^2)^3}$ & $\frac{1}{(1+\alpha^2u^2)^3}$ & $\frac{s^5}{(s^2+\alpha^2)^3}$ \\
  28 & $\frac{t^2\sin(\alpha t)}{2\alpha}$ & $\frac{u^4(3s^2-\alpha^2u^2)}{(s^2+\alpha^2u^2)^3}$ & $\frac{u^3(3s^2-\alpha^2u^3)}{(s^2+\alpha^2u^2)^3}$ & $\frac{u^3(-3+\alpha^2u^2)}{(1+\alpha^2u^2)^3}$ & $\frac{(3s^2-\alpha^2)}{(s^2+\alpha^2)^3}$ \\
  29 & $\frac{t^2\cos(\alpha t)}{2}$ & $\frac{u^3(s^3-3\alpha^2u^2s)}{(s^2+\alpha^2u^2)^3}$ & $\frac{u^2(s^3-3\alpha^2u^2s)}{(s^2+\alpha^2u^2)^3}$ & $\frac{u^2(1-3\alpha^2u^2)}{(1+\alpha^2u^2)^3}$ & $\frac{(s^3-3\alpha^2s)}{(s^2+\alpha^2)^3}$ \\
  30 & $\frac{t^3\sin(\alpha t)}{24\alpha}$ & $\frac{su^5(s-\alpha u)^2}{(s^2+\alpha^2u^2)^4}$ & $\frac{su^4(s-\alpha u)^2}{(s^2+\alpha^2u^2)^4}$  & $\frac{u^4(1-\alpha u)^2}{(1+\alpha^2u^2)^4}$  & $\frac{s(s-\alpha)^2}{(s^2+\alpha^2)^4}$  \\
  31 & $\frac{\exp(\alpha t)-\exp(\beta t)}{\alpha-\beta}\,\,\alpha\neq \beta$ & $\frac{u^2}{(s-\alpha u)(s-\beta u)}$ & $\frac{u}{(s-\alpha u)(s-\beta u)}$ & $\frac{u}{(1-\beta u)(1-\alpha u)}$ & $\frac{1}{(s-\beta)(s-\alpha)}$ \\
  32 & $\frac{\alpha \exp(\alpha t)-\beta \exp(\beta t)}{\alpha-\beta}\,\,\alpha\neq \beta$ & $\frac{us}{(s-\beta u)(s-\alpha u)}$ & $\frac{s}{(s-\beta u)(s-\alpha u)}$ & $\frac{1}{(1-\beta u)(1-\alpha u)}$ & $\frac{s}{(s-\beta)(s-\alpha)}$ \\
  33 & $I_0(\alpha t)$ & $\frac{u}{\sqrt{s^2-\alpha^2 u^2}}$ & $\frac{1}{\sqrt{s^2-\alpha^2 u^2}}$ & $\frac{1}{\sqrt{1-\alpha^2 u^2}}$ & $\frac{1}{\sqrt{s^2-\alpha^2}}$ \\
  34 & $\delta(t-\alpha)$ & $u\exp\left(\frac{-\alpha s}{u}\right)$ & $\frac{1}{u}\exp\left(\frac{-\alpha s}{u}\right)$ & $\frac{1}{u}\exp\left(\frac{-\alpha}{u}\right)$ & $\exp(-\alpha s)$ \\
  35 & $J_0(\alpha t)$ & $\frac{u}{\sqrt{s^2+\alpha^2 u^2}}$ & $\frac{1}{\sqrt{s^2+\alpha^2 u^2}}$ & $\frac{1}{\sqrt{1+\alpha^2 u^2}}$ & $\frac{1}{\sqrt{s^2+\alpha^2}}$ \\
  \hline
\end{longtable}}
\end{center}
\vspace{1cm} \textbf{Table 2:} General properties of Shehu
transform {\setlength{\extrarowheight}{10pt}
\begin{longtable}{|c|c|c|c|}
  \hline
  S.NO. & Definition/Property & ${\mathbb S}\left[v(t)\right]=V(s,u)$ & Transforms \\
  \hline
  1 & Definition & ${\mathbb S}\left[v(t)\right]$ & $\int _{0 }^{\infty
}\exp\left(\frac{-st}{u}\right)v(t)\, dt;\, \, \, \, s\,>0,\,\, u>0$ \\
  2 & Inverse & $v(t)={\mathbb S^{-1}}\left[V(s,u)\right]$ & $\frac{1}{2\pi i}\int
_{\alpha-i\infty}^{\alpha+i\infty}\frac{1}{u}\exp\left(\frac{st}{u}\right)V(s,u)\, ds$ \\
  3 & Linearity & ${\mathbb S}\left[\alpha v(t)+\beta w(t)\right]$ & $\alpha{\mathbb S}\left[v(t)\right]+\beta{\mathbb S}\left[w(t)\right]$ \\
  4 & Change of scale & ${\mathbb S}\left[v(\alpha t)\right]$ & $\frac{u}{\alpha}V\left(\frac{s}{\alpha},u\right)$ \\
  5 & Derivatives & ${\mathbb S} \left[v'(t)\right]$ & $\frac{s}{u} V(s,u)-v(0)$ \\&  & ${\mathbb S}\left[v''(t)\right]$ & $\frac{s^{2} }{u^{2} } V(s,u)-\frac{s}{u}v(0)-v'(0)$ \\ &  & ${\mathbb S}\left[v'''(t)\right]$ & $\frac{s^{3} }{u^{3} }V(s,u)-\frac{s^2}{u^2}v(0)-\frac{s}{u}v'(0)-v''(0)$ \\ &  & $\vdots$ &  \\ &  & ${\mathbb S} \left[v^{(n)}(t)\right]$ & $\frac{s^n}{u^n}V(s,u)-\sum_{k=0}^{n-1}\left(\frac{s}{u}\right)^{n-(k+1)}$\\
  \hline
\end{longtable}}
\vspace{8cm} \textbf{Table 3:} Summary of some integral transform
and their definitions {\setlength{\extrarowheight}{10pt}
\begin{longtable}{|c|c|c|}
  \hline
  S.No. & Integral transform & Definition \\
  \hline
  1 & Laplace transform & $\pounds[f(t)]=F(s)=\int_{0}^{\infty}\exp\left(-st\right)f(t)dt$ \\
  2 & Fourier  transform & $\digamma[f(t)]=f(\omega)=\frac{1}{\sqrt{2\pi}}\int_{-\infty}^{\infty}\exp\left(-i\omega t\right)f(t)dt$ \\
  3 & Melling transform & $M[f(s);s]=f^{*}(s)=\int_{0}^{\infty}x^{s-1}f(x)dx$ \\
  4 & Hankel's transform & $F_v(s)=H_v[f(r)]=\int_{0}^{\infty}rf(r)J_v(sr)dr,\,\,\,\,\,r\geq 0$ \\
  5 & Sumudu transform & $\textbf{S}[f(t)]=G(u)=\frac{1}{u}\int_{0}^{\infty}\exp\left(\frac{-t}{u}\right)f(t)dt$ \\
  6 & Laplace-Carson transform & $\hat{f}_C(p)=p\int_{0}^{\infty}\exp\left(-pt\right)f(t)dt$,\,\,$t\geq0$ \\
  7 & Atangana-Kilicman transform & $M_n(s)=M_n[f(x)](s)=\int_{0}^{\infty}x^n\exp\left(-xs\right)f(x)dx$\\
  8 & El-zaki transform & $E[f(t)]=T(u)=u\int_{0}^{\infty}\exp\left(\frac{-t}{u}\right)f(t)dt$ \\
  9 & Yang transform & $Y[\phi(\tau)]=\phi(\omega)=\int_{0}^{\infty}\exp\left(\frac{-\tau}{\omega}\right)\phi(\tau)d\tau$ \\
  10 & natural transform & ${\mathbb N^+}[f(t)]=R(s,u)=\frac{1}{u}\int_{0}^{\infty}\exp\left(\frac{-st}{u}\right)f(t)dt,\, \, \, \, s\,>0,\,\, u>0$ \\
  11 & z-transform & $X(z)=Z\{x[n]\}=\sum_{n=-\infty}^{\infty}x[n]z^{-n},\,\,n\in{\mathbb Z},\,\,z\in{\mathbb C}$\\
  12 & ${\mathbb M}$-transform & $M_{\rho,m}[f(t)](u,v)=\int_{0}^{\infty}\frac{\exp\left(-ut\right)f(vt)}{\left(t^m+v^m\right)^\rho}dt,\,\,\rho \in{\mathbb C},\,\Re(\rho)\geq0,\,m\in {\mathbb
Z_+}$ \\
  \hline
\end{longtable}}


\begin{thebibliography}{10}
%
\baselineskip=14pt
\bibitem{R1}
H.A. Agwa, F.M. Ali, A. Kilicman,
\newblock A new integral transform on time scales and its applications,
\newblock {\it Advances in Difference Equations}, {\bf 60}(2012), 1--14.


\bibitem{R2}
C. Ahrendt,
\newblock The Laplace transform on time scales,
\newblock {\it Pan. Am. Math. J.}, {\bf 19}(2009), 1--36.

\bibitem{R3}
A. Atangana,
\newblock A note on the triple Laplace transform and its applications to some kind of third-order differential equation,
\newblock {\it Abstract and Applied Analysis}, {\bf 2013}(2013), Article ID 769102, 1--10.

\bibitem{R4}
H.M. Srivastava, A.K. Golmankhaneh, D. Baleanu, X.Y. Yang,
\newblock Local
fractional Sumudu transform with applications to IVPs on Cantor
sets,
\newblock {\it Abstr Appl Anal}, {\bf 2014}(2014), Article ID
620529,
1--7.

\bibitem{R5}
G. Dattoli, M. R. Martinelli, P. E. Ricci,
\newblock On new families of integral transforms for the solution of partial differential equations,
\newblock {\it Integral Transforms and Special Functions}, {\bf 8}(2005),
661--667.

\bibitem{R6}
H. Bulut, H.M. Baskonus, and F.B.M. Belgacem,
\newblock The analytical solution of some fractional ordinary differential equations by the Sumudu transform method,
\newblock {\it Abstract and Applied Analysis}, {\bf 2013}(2013), Article ID 203875, 1--6.

\bibitem{R7}
S.Weerakoon,
\newblock The Sumudu transform and the Laplace transform: reply,
\newblock {\it International Journal of Mathematical Education in Science
and Technology}. {\bf 28}(1997), 159--160.

\bibitem{R8}
D. Albayrak, S.D. Purohit, and F. U\c{c}ar,
\newblock Certain inversion and representation formulas for q-sumudu transforms,
\newblock {\it Hacettepe Journal of Mathematics and Statistics}, {\bf 43}(2014), 699--713.

\bibitem{R9}
S. Weerakoon,
\newblock Application of Sumudu transform to partial
differential equations,
\newblock {\it International Journal of Mathematical
Education in Science and Technology,} {\bf 25}(1994), 277--283.

\bibitem{R10}
X.Y. Yang, Y. yang, C. Cattani, and M. Zhu,
\newblock A new technique for solving 1-D Burgers equation,
\newblock {\it Thermal Science,} {\bf 21}(2017), S129--S136.

\bibitem{R11}
A. Kilicman, H. Eltayeb,
\newblock On new integral transform and differential equations,
\newblock {\it J. Math. Probl. Eng,} {\bf 2010}(2010), Article ID: 463579, 1--13.

\bibitem{R12}
S. Bochner, K. Chandrasekharan,
\newblock Fourier transforms,
\newblock {\it Princeton University Press, Princeton, NJ, USA}, (1949).

\bibitem{R13}
R. N. Bracewell,
\newblock The Fourier transform and its applications,
\newblock {\it McGraw-Hill, Boston, Mass, USA, 3rd edition}, (2000).

\bibitem{R14}
R. Murray, Spiegel.
\newblock Theory and problems of Laplace transform.
\newblock {\sl New York, USA: Schaum's Outline Series, McGraw--Hill,} (1965).

\bibitem{R15}
L. Debnath, D. Bhatta..
\newblock Integral transform and their applications.
\newblock {\sl CRC Press, New York, NY, USA} (2010).

\bibitem{R16}
B. Davies,
\newblock Integral transforms and their applications,
\newblock {\it Texts in Applied Mathematics, Springer, New York, NY, USA}, {\bf}(2002).

\bibitem{R17}
E.I. Jury,
\newblock Theory and applications of the z-transform
Method,
\newblock {\it John Wiley and Sons, New York, NY,USA,} (1964).

\bibitem{R18}
K. Liu, R.J. Hu, C. Cattani, G.N. Xie, X.J. Yang, and Y. Zhao,
\newblock Local fractional z-transforms with applications to
signals on Cantor sets,
\newblock {\it Abstract and Applied Analysis}, {\bf 2014}(2014),Article ID: 638648, 1--6.

\bibitem{R19}
 P.M. Morse, H. Feshbach,
\newblock Methods of theoretical physics,
\newblock {\it McGraw-Hill, New York}, (1953), 484--485.


\bibitem{R20}
P. Flajolet, X. Gourdon, and P. Dumas,
\newblock Mellin transforms and asymptotics: harmonic sums,
\newblock {\it Theoretical Computer Science}, {\bf 144}(1995), 3--58.

\bibitem{R21}
C. Donolato,
\newblock Analytical and numerical inversion of the Laplace-Carson
transform by a differential method,
\newblock {\it Computer Physics Communications,} {\bf 145}(2002), 298--309.

\bibitem{R22}
A.M. Makarov,
\newblock Application of the Laplace-Carson method of integral transformation to the theory of unsteady visco-plastic flows,
\newblock {\it J. Engrg. Phys. Thermophys} {\bf 19}(1970), 94--99.

\bibitem{R23}
E. Sjøntoft,
\newblock A straightforward deconvolution method for use in small computers,
\newblock {\it Nucl. Instrum. Methods}, {\bf 163}(1979), 519--522.


\bibitem{R24}
A.S. Vasudeva Murthy,
\newblock A note on the differential inversion method of
Hohlfield et al.,
\newblock {\it SIAM J. Appl. Math.}, {\bf 55}(1995), 712--722.

\bibitem{R25}
I.N. Sneddon,
\newblock The Use of integral transform,
\newblock {\it McGraw-Hill, New York}, (1972).

\bibitem{R26}
K. Xie, Y. Wanga, K. Wang, and X. Cai,
\newblock Application of Hankel transforms to boundary value problems of water
flow due to a circular source,
\newblock {\it Applied Mathematics and Computation}, {\bf 216}(2010), 1469--1477.

\bibitem{R27}
Watugala GK.
\newblock Sumudu transform--a new integral transform to solve differential equations and control engineering
problems.
\newblock  {\sl Math. Engg. in Indust.,} {\bf 6}(1998), 319--329.

\bibitem{R28}
M.A. Asiru,
\newblock  Sumudu transform and solution of integral equations of convolution type,
\newblock{\sl Int. J. of Math. Edu. Sci. and Tech.,} {\bf 33}(2002), 944--949.

\bibitem{R29}
 F.B.M. Belgacem, S.L. Kalla, A.A. Karaballi,
\newblock Analytical investigations of the Sumudu transform and applications to integral production
equations,
\newblock{\sl Math. Probl. in Engg.}, {\bf 3}(2003),
103-118.

\bibitem{R30}
 F.B.M. Belgacem, A.A. Karaballi,
\newblock Sumudu transform fundamental properties, investigations and
applications.
\newblock {\sl J. of Appl. Math. and Stoch. Anal.,} {\bf 2006}(2006), Article ID 91083, 1--23.

\bibitem{R31}
H. ELtayeh and A. kilicman,
\newblock On Some applications of a new integral
transform,
\newblock {\it Int. Journal of Math. Analysis}, {\bf 4}(2010), 123--132.

\bibitem{R32}
T.M. Elzaki.
\newblock The new integral transform ''Elzaki transform''.
\newblock {\sl Glob. J. of Pur. and Appl. Math.,} {\bf 7}(2011), 57--64.

\bibitem{R33}
Z.H. Khan, W.A. Khan,
\newblock N-transform-properties and applications.
\newblock {\sl NUST J. of Engg. Sci.,} {\bf 1}(2008), 127--133.

\bibitem{R34}
F.B.M. Belgacem, R. Silambarasan,
\newblock  Theory of natural transform.
\newblock{\sl Math. in Engg. Sci., and Aeros.,} {\bf 3}(2012), 99--124.

\bibitem{R35}
F.B.M. Belgacem, R. Silambarasan,
\newblock Advances in the natural transform.
\newblock {\sl AIP Conference Proceedings; 1493 January 2012; USA: American Institute of Physics.} (2012), 106--110.

\bibitem{R36}
H.M. Srivastava, Minjie Luo, R.K.Raina,
\newblock A new integral transform and its applications,
\newblock {\it Acta Mathematica Scientia}, {\bf 35}(2015), 1386--1400.

\bibitem{R37}
A. Atangana, A. Kilicman,
\newblock A novel integral operator transform and its application to
some FODE and FPDE with some kind of singularities,
\newblock {\it Mathematical Problems in Engineering}, {\bf 2013}(2013), Article ID: 531984, 1--7.

\bibitem{R38}
X.J. Yang,
\newblock A new integral transform method for solving steady heat-transfer problem,
\newblock {\it Thermal Science}, {\bf 20}(2016), S639--S642.

\bibitem{R39}
X.J. Yang,
\newblock A new integral transform operator for solving the heat-diffusion problem,
\newblock {\it Applied Mathematics Letters}, {\bf 64}(2017), 193--197.

\bibitem{R40}
Y.X. Jun, F. Gao ,
\newblock A new technology for solving diffusion and heat equations,
\newblock {\it Thermal Science}, {\bf 21}(2017), 133--140.

\bibitem{R41}
A. Atangana, B.S.T. Alkaltani,
\newblock A novel double integral transform and its
applications,
\newblock {\it Journal of Nonlinear Science and Applications}, {\bf 9}(2016), 424--434.

\bibitem{R42}
H. Eltayeb,
\newblock A note on double Laplace decomposition method and
nonlinear partial differential equations,
\newblock {\it New Trends in Mathematical Sciences}, {\bf 5}(2017), 156--164.

\bibitem{R43}
F.B.M. Belgacem, R. Silambarasan, H. Zakia, T. Mekkaoui,
\newblock New and extended applications of the natural and Sumudu transforms: Fractional diffusion and Stokes fluid flow realms.
\newblock {\it Advances Real and Complex Analysis with Applications, Publisher: Birkhäuser, Singapore,} (2017).
107--120.
\end{thebibliography}
\end{document}